\documentclass[10pt]{amsart}

\usepackage{amsmath}

\usepackage[lite]{amsrefs}

\usepackage{graphicx}

\usepackage{amssymb}

\usepackage[all,cmtip]{xy}

\newcommand{\R}{\mathbb{R}}

\newcommand{\Z}{\mathbb{Z}}

\newcommand{\vcenteredinclude}[1]{\begingroup
\setbox0=\hbox{\includegraphics[scale=0.6]{#1}}%
\parbox{\wd0}{\box0}\endgroup}

\newcommand{\vcenteredincludeb}[1]{\begingroup
\setbox0=\hbox{\includegraphics[scale=0.4]{#1}}%
\parbox{\wd0}{\box0}\endgroup}

\newcommand{\vcenteredincluded}[1]{\begingroup
\setbox0=\hbox{\includegraphics[scale=0.4]{#1}}%
\parbox{\wd0}{\box0}\endgroup}

\newcommand{\vcenteredincludee}[1]{\begingroup
\setbox0=\hbox{\includegraphics[scale=0.5]{#1}}%
\parbox{\wd0}{\box0}\endgroup}

\numberwithin{equation}{section}

\theoremstyle{plain} 

\newtheorem{thm}[equation]{Theorem}

\newtheorem{prop}[equation]{Proposition}

\theoremstyle{definition}

\newtheorem{defn}[equation]{Definition}

\theoremstyle{remark}

\newtheorem{ex}[equation]{Example}

\usepackage{mathtools}

\usepackage{graphicx}

\usepackage{subcaption}

\begin{document}






\title{Final Report\\Knots, Braids, and Knot Invariants}






\author{Matthew Stevens \\ Consultant: Professor Matt Hogancamp \\ Fall 2021} \email{stevens.mat@northeastern.edu}










\date{\today}



\begin{abstract}

In this report, I will start by first giving a brief introduction on knots to build some intuition before beginning the more rigorous review in the Literature Review section. There, I will define knot equivalence, the Jones polynomial invariant, braid groups, Hecke algebras, and the HOMFLY polynomial invariant. Next, I will cover some of my independent research regarding the computation of a general formula for the HOMFLY polynomial of the closure of looped coxeter braids. Lastly, I will loosely discuss how this research connects to future research regarding the universal trace in the Plan for Future Work section.

\end{abstract}

\maketitle


\tableofcontents

\pagebreak




\section{Introduction}

\label{sec:intro}

I will define this more rigorously below, but a knot is usually referred to as a closed loop in 3-dimensional space, while a link would be a collection of closed loops in 3-dimensional space. The only important stipulation is that these loops never intersect themselves or other closed loops. Figure 1 below displays some notable knots and links:

\begin{figure}[h]
    \centering
    \includegraphics[width=.8\textwidth]{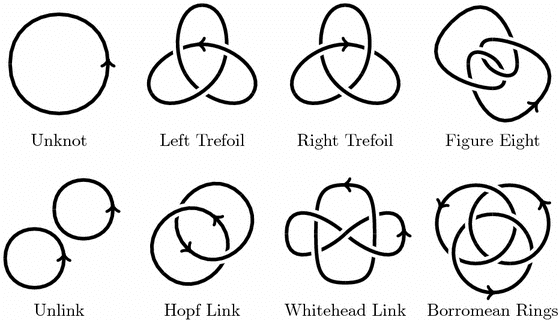}
    \caption{Examples of Knots and Links [5]}
\end{figure}

As seen, knots and links can be projected into a two dimensional space given a convention for displaying when a section of the loop crosses over or under another section of the loop. The arrows in the above picture determine the orientation of the loop. The orientation of a knot or link will be useful in defining knot invariants. In 1927, Kurt Reidemeister proved that two knot diagrams belonging to the same knot can be related in some sequence of three moves shown in Figure 2 below: 

\begin{figure}[h]
    \centering
    \includegraphics[width=.4\textwidth]{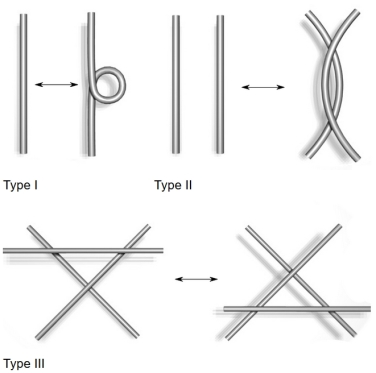}
    \caption{The Three Reidemeister Moves [6]}
\end{figure}

Two knots are equivalent if you can continuously change one into the other without intersecting any sections of the loops during that change. This notion of knot equivalence will be defined more strictly in the next section. Figure 3 displays four knots that are all equivalent to each-other.

\begin{figure}[h]
    \centering
    \includegraphics[width=.8\textwidth]{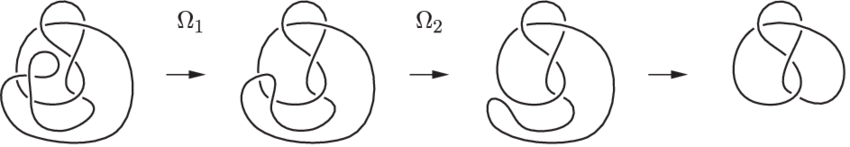}
    \caption{Equivalent Knots [7]}
\end{figure}

As you can see, the 2-dimensional representations of these knots are slightly different. For instance, the knot on the left contains three more crossings than the knot on the far right. Due to this property of equivalent knots, there has been much work in mathematics to discover knot invariants. These knot invariants are properties of a knot that does not change when the knot is manipulated by any of the three Reidemeister moves shown above. These invariants provide a simple way to determine if two knots are equivalent. It will be seen in the below section that the Jones polynomial and HOMFLY polynomial are types of knot invariants.






















\section{Literature Review}

\subsection{Links and Knots}

In this section, I will cover more rigorous definitions of knots and knot equivalence.

\begin{defn}

For any two distinct points, \(p,q\in\R^3 \), let \([p,q]\) denote the linear line segment joining them. Then, for an ordered set of distinct points, \((p_1,p_2,...,p_n)\in\R^3\), we have that the union of the linear line segments formed from these points, \([p_1,p_2],[p_2,p_3],...,[p_{n-1},p_n],[p_n,p_1]\), is a closed polygonal curve. If each line segment only intersects other line segments in two places, each at the two endpoints of the segment, then the curve is said to be simple.

\end{defn}

\begin{defn}

A Link \(L\) of \(m\in\mathbb{N}\cup\{0\}\) components is a subset of \(\R^3\) that consists of \(m\) disjoint, polygonal, simple closed curves. A knot is a link with one component. A knot will be determined by a single sequence as seen in Definition 2.1, while a link will be determined by a set of these sequences.

\end{defn}

As seen in Definition 2.1, polygonal means that the curves are composed of a finite number of straight lines. The finite number of lines is important so as to avoid complications involving infinite kinks that converge to a single point. It is usually assumed that there is so many straight lines composing a loop that it appears smooth when drawn. This can be seen in the picture of various knots and links provided in the introduction. Figure 4 provides imagery of what polygonal knots, composed of a small number of line segments, look like. The unknot, as seen in Figure 4a, can be composed of a minimum of three line segments.

\begin{figure}

\centering
\begin{subfigure}{0.4\textwidth}
    \includegraphics[width=0.8\textwidth]{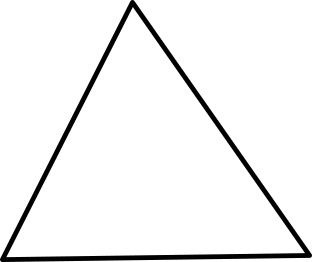}
    \caption{The unkot}
    
\end{subfigure}
\hfill
\begin{subfigure}{0.4\textwidth}
    \includegraphics[width=0.8\textwidth]{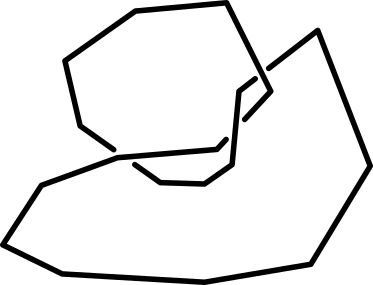}
    \caption{The trefoil knot}
    
\end{subfigure}
\hfill
        
\caption{Knot examples from Definition 2.2}

\end{figure}

\begin{defn}

A link \(J\) is an elementary deformation of a link \(K\) if one link is determined by a set of sequences \(\{(p_{1_{1}},p_{1_{2}},...p_{1_{n}}),...,(p_{m_{1}},p_{m_{2}},...p_{m_{l}}) \}\), and the other link is determined by the sequence \(\{(p_{1_{0}},p_{1_{1}},p_{1_{2}},...p_{1_{n}}),...,(p_{m_{1}},p_{m_{2}},...p_{m_{l}}) \}\) given the following conditions:

(i) \(p_{1_{0}}\) is not collinear with \(p_{1_{1}}\) and \(p_{1_{n}}\)

(ii) The triangle spanned by \((p_{1_{0}},p_{1_{1}},p_{1_{n}})\) only intersects the link determined by the sequences \(\{(p_{1_{1}},p_{1_{2}},...p_{1_{n}}),...,(p_{m_{1}},p_{m_{2}},...p_{m_{l}}) \}\)  in the line segment \([p_{1_{n}},p_{1_{1}}]\)

\end{defn}

Figure 2 displays a link deformation that would not be allowed. As seen, the triangle spanned by adding the new point intersects a segment of the link. This deformation takes the deformed segment from being over to being under in regards to the crossing.

\begin{figure}

    \centering
    \includegraphics[width=0.8\textwidth]{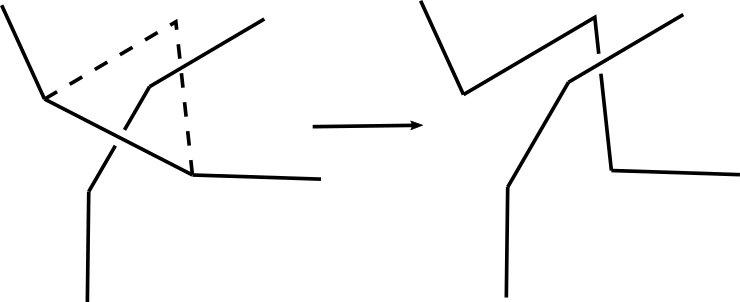}
    \caption{Example of an unwanted deformation}
    
\end{figure}

\begin{defn}

Two links, K and J, are called equivalent if there exists a sequence of links \((K_0,K_1,...,K_n)\) such that \(K_0=K\), \(K_n=J\), and the link \(K_{i+1}\) is an elementary deformation of the link \(K_i\) for \(i\in\{0,1,2,..,n-1\}\)

\end{defn}

As referenced in the introduction, it is useful to consider two dimensional representations of links as these are simpler to draw. In fact, the majority of the links featured have been two dimensional representations of three dimensional links. The following definition describes the most common way to project a link from \(\R^3\) into \(\R^2\). 

\begin{defn}

Let K be a link in \(\R^3\). The projection of K is the image of K under the function \(\pi: \R^3 \to \R^2\) such that \(\pi(x,y,z)=(x,y)\).

\end{defn}

The projection of a link K will look like a set of closed curves in \(\R^2\) that may cross itself. The crossings in the projection indicate when one segment of the link crosses over or under another segment of the link. In order to illustrate this lost information in the link projection, a gap in one of the segments during a crossing is used to show which of the two segments crosses under the other segment. The projection of a link combined with this extra information at the crossings is called a link diagram (or knot diagrams when dealing with single component links). The Figures 4a and 4b are knot diagrams of the unknot and trefoil knot.

\begin{defn}

The link diagram of a link \(K\) is called a regular diagram if no three points of the link project to the same point on the diagram and if no vertex point on the link projects to the same point as another point on the link.

\end{defn}

The point of this specification is to avoid link diagrams that have limited information due to specific projections. Some of these unwanted projections can be seen in Figure 6. As seen, having more than three segments crossing at the same point makes it difficult to distinguish the crossing heights of the segments. Moreover, allowing vertices to project to the same point as other points on a link creates a few issues. These include projecting two different segments to the same points in \(\R^2\). Going forward, regular link/knot diagrams will be the only link/knot diagrams of importance. 

\begin{figure}
    \centering
    \includegraphics[width=.8\textwidth]{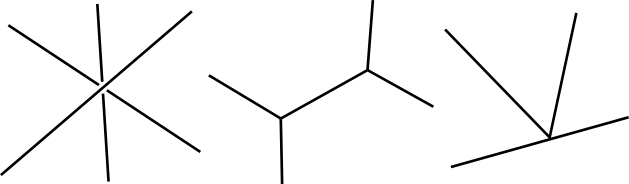}
    \caption{Three examples of unwanted projections}
\end{figure}

\begin{defn}

An oriented knot consists of a knot and an ordering of its vertices. An oriented link is a union of disjoint oriented knots.

\end{defn}

As an example, the orientation of an oriented knot \(K\) is determined by the ordering of its vertices \((p_1,p_2,...,p_n)\). Any cyclic permutation of these vertices, such as \((p_n,p_1,...,p_{n-1})\), produces an equivalent oriented knot. In contrast, the reversal of the ordering of the vertices, such as \((p_n,...,p_2,p_1)\), will produce the same knot but the oriented knot is not necessarily equivalent. The knot diagram of an oriented knot will be drawn with an arrow representing the orientation which can be seen in Figure 7.  Oriented knots and links will become useful when discussing knot invariants. 

\begin{figure}
    \centering
    \includegraphics[width=0.4\textwidth]{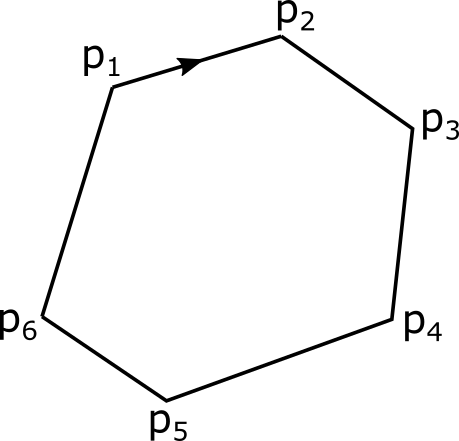}
    \caption{Oriented unknot}
\end{figure}

When considering a knot K, an elementary deformation of this knot can drastically alter its knot diagram. Because of this, it is helpful to break down the most basic ways a knot diagram can be altered using elementary deformations. Figure 8 illustrates three of these basic changes in knot diagrams that can occur through elementary deformations of a knot. These three moves are known as Reidemeister moves.

\begin{thm}[Reidemeister's Theorem]

If two knots, \(K\) and \(J\), are equivalent and have regular diagrams, then their knot diagrams can be related by a sequence of Reidemeister moves. Likewise, if two knots, \(K\) and \(J\), have regular diagrams and their diagrams can be related by sequence of Reidemeister moves, then these knots are equivalent. This theorem can also be extended to links.

\end{thm}

Theorem 2.8 will be extremely powerful in helping to determine invariants of knots. If a certain property of a knot does not change under the three Reidemeister moves, then this property is an invariant of the knot. In this case, invariant means that the property is constant for all knots in an equivalence class defined under Definition 2.4.

\begin{figure}
    \centering
    \includegraphics[width=0.6\textwidth]{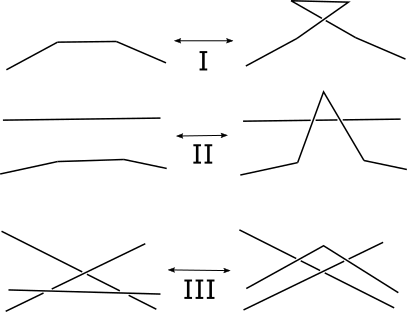}
    \caption{Reidemeister moves}
\end{figure}

\subsection{The Jones Polynomial}

 A Jones polynomial can be calculated in a recursive way given any link diagram. This is done using skein relations ((i), (ii), and (iii) in Definition 2.9) that allow one to reduce links down to links in which the Jones polynomial is known.

\begin{defn}

The Jones polynomial invariant is a function 

\[V: \{\text{Oriented links in } \R^3\} \to \Z[t^{-1/2},t^{1/2}]\] 

satisfying the following skein relations:

(i) If two knots are equivalent then their Jones polynomials are the same.

(ii) V(unknot) = 1

(iii) Whenever three oriented links, \(L_+\), \(L_-\), and \(L_0\), differ only in the neighborhood of a point as shown in the picture below, then

\[t^{-1}V(L_+) - tV(L_-) + (t^{-1/2} - t^{1/2})V(L_0)=0\]

\begin{center}
    \includegraphics[width=.6\textwidth]{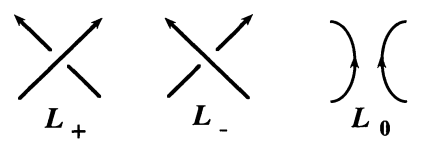}
\end{center}

\end{defn}

\begin{prop}

The Jones polynomial of the split union of a link L and the unknot is

\[V(L\sqcup unknot) = \frac{t-t^{-1}}{t^{-1/2}-t^{1/2}}V(L)\]

\end{prop}

Here the split union refers to a link composed of the link L and the unknot, where the link diagram of \(L\sqcup unknot\) can be drawn so that there is no crossing points between the unknot and the link L.

\bigskip

Proof: We can consider the link diagram of the split union of a general link L with the unknot. The link L is represented by a dotted box and three segments jutting out of the box. This is all that is needed to understand why Proposition 2.10 is true. We have that

\[t^{-1}V(\includegraphics[height=0.13\textwidth]{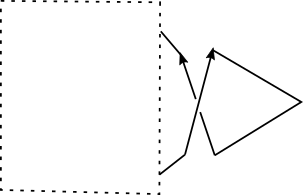})-tV(\includegraphics[height=0.13\textwidth]{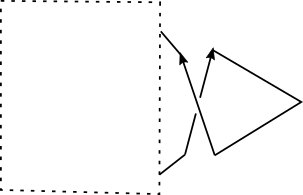})+(t^{-1/2}-t^{1/2})V(\includegraphics[height=0.13\textwidth]{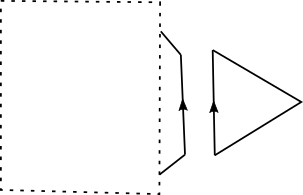})=0.\]

We can see that the first two diagrams are equivalent to the link L through the first Reidemeister move. Therefore, we have that

\[t^{-1}V(L)-tV(L)+(t^{-1/2}-t^{1/2})V(L\sqcup unknot)\].

With some manipulation we get the equation seen in Proposition 2.10. \(\square\)

\begin{ex}

This example will cover the Jones polynomial calculation of an oriented hopf link:

\[t^{-1}V(\includegraphics[height=0.1\textwidth]{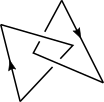})-tV(\includegraphics[height=0.1\textwidth]{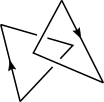})+(t^{-1/2}-t^{1/2})V(\includegraphics[height=0.1\textwidth]{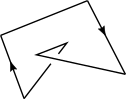})=0\]

Due to Theorem 2.8 and Definition 2.9, we know that

\[V(\includegraphics[height=0.1\textwidth]{hopf3a.png})=1.\]

We also know that

\[V(\includegraphics[height=0.1\textwidth]{hopf2a})=V(\includegraphics[height=0.1\textwidth]{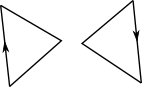}).\]

Now, through are definition of equivalence and Definition 2.9, we have that

\[V(\includegraphics[height=0.1\textwidth]{2unkot.png})=V(\includegraphics[height=0.1\textwidth]{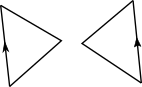}).\]

Moreover, through proposition 2.10, we have that

\[V(\includegraphics[height=0.1\textwidth]{2unkot2.png})=\frac{t-t^{-1}}{t^{-1/2}-t^{1/2}}V(unknot).\]

Putting it all together, we get that

\[t^{-1}V(\includegraphics[height=0.1\textwidth]{hopfa.png})-t\frac{t-t^{-1}}{t^{-1/2}-t^{1/2}}+(t^{-1/2}-t^{1/2})=0.\]

Finally, we get that

\[V(\includegraphics[height=0.1\textwidth]{hopfa.png})= -t^{5/2}-t^{1/2}.\]

\end{ex}

As seen in the example above, you can use the skein relation (iii) to relate the Jones polynomial of the desired link to link diagrams of split unions of the unknot. From there, you can use the result from Proposition 2.10 in order to relate the Jones polynomial of the desired link to a link diagram of the unknot. Lastly, you can use the skein relation (ii) in order to complete the calculation of the Jones polynomial of the desired link.

\subsection{Braid Groups}

A mathematical braid of rank \(n \in \mathbb{N}\) is composed of \(n\) strands that permute \(n\) nodes. Figure 9 displays an example of a mathematical braid of rank 4. All braids of rank 4 form a group under the operation of braid concatenation. 

\begin{figure}[h]
    \centering
    \includegraphics[width=0.3\textwidth]{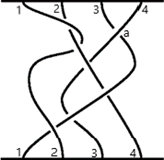}
    \caption{Braid Example}
\end{figure}

\begin{figure}[h]
    \centering
    \includegraphics[width=0.3\textwidth]{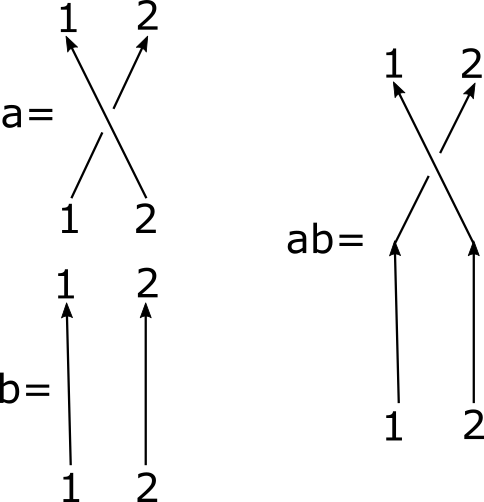}
    \caption{Braid Concatenation}
\end{figure}

Figure 10 displays an example of what braid concatenation looks like for braids of rank 2. As seen, you connect the top nodes of b to the corresponding bottom nodes of a, forming another braid of rank 2. It should be clear that the set of braids of a specific rank is closed under this operation. An algebraic representation of the braid group of a specific rank is given in Definition 2.12.

\begin{defn}

\(B_n = \langle \sigma_1,...,\sigma_{n-1} | \sigma_i\sigma_j = \sigma_j\sigma_i \forall i,j s.t. |i-j| \geq 2; \sigma_i\sigma_j\sigma_i = \sigma_j\sigma_i\sigma_j \forall i,j s.t. |i-j|=1\rangle \)

\end{defn}

\begin{figure}[h]
    \centering
    \includegraphics[width=0.5\textwidth]{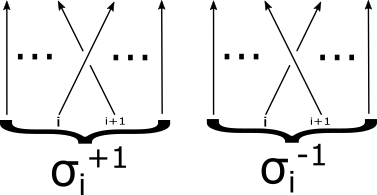}
    \caption{Braid Representation of Generators}
\end{figure}

\begin{figure}[h]
    \centering
    \includegraphics[width=0.4\textwidth]{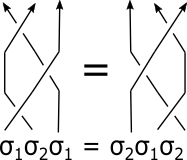}
    \caption{Braid Representation of Second Relation}
\end{figure}

As seen from Definition 2.12, the braid group of rank n can be generated by n-1 generators while enforcing two relations. The generator \(\sigma_i\), displayed in Figure 11, can be thought of as the braid of rank n where only one crossing occurs between the strand starting at the node i and the strand starting at the node i+1. It's also important to note that this crossing is positive, meaning that the strand starting on the left crosses over the strand starting at the right. The same braid as \(\sigma_i\) except with a negative crossing as opposed to a positive crossing is actually the inverse of \(\sigma_i\). Lastly, Figure 12 illustrates the second relation of our braid group presentation on a braid of rank 3. Take notice of the similarities between this relation and the type III Reidemeister move. 

Visually, it should be clear that their seems to be obvious similarities between braids and knots/links. In fact, there exists an operation referred to as a braid closure which allows one to produce a link from a braid.

\begin{defn}

The closure of a braid occurs by connecting each corresponding nodes in pairs. It is important that no extra crossing is produced through this closure. Figure 13 displays a general braid closure while Figure 14 displays an example braid closure on a rank 3 braid.

\end{defn}

\begin{figure}[h]
    \centering
    \includegraphics[width=0.5\textwidth]{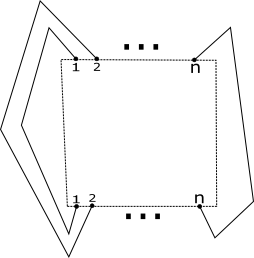}
    \caption{General Braid Closure}
\end{figure}

\begin{figure}[h]
    \centering
    \includegraphics[width=0.4\textwidth]{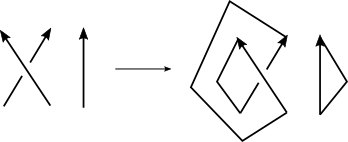}
    \caption{Braid Closure Example}
\end{figure}

Theorem 2.14, proved by James Waddell Alexander II in 1923, explicitly states the connection between braids and links. This theorem will become important in the next section regarding Hecke algebras.

\begin{thm}

Every knot or link can be represented as the closure of a braid.

\end{thm}

\subsection{Hecke Algebras}

\begin{defn}

\(\mathbb{H}_n = \mathbb{Z}[q,q^{-1}]\)-algebra generated by n-strand braids modulo the following relations:

1. \(\vcenteredinclude{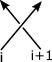} - \vcenteredinclude{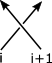} = (q-q^{-1})\vcenteredinclude{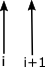}\)

2. \(\vcenteredincludeb{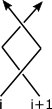} = \vcenteredinclude{1c.png}\)

3. \(\vcenteredincludeb{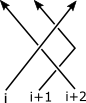} = \vcenteredincludeb{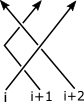}\)

\end{defn}

As seen in Definition 2.15, the Hecke Algebra of rank n can be thought of as a R-Algebra with the ring being the ring of Laurent polynomials, \(\mathbb{Z}[q,q^{-1}]\). The multiplication defined on this algebra is as follows:

\[(\sum_{i=0}^{m} f_i \vcenteredincludee{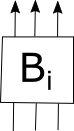})(\sum_{j=0}^{n} f_j \vcenteredincludee{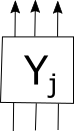})= \sum_{i=0}^{m}\sum_{j=0}^{n}f_if_j\vcenteredincludee{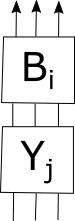}, f_i,f_j \in \mathbb{Z}[q,q^{-1}], B_i,Y_j \in B_n\]

As seen, the elements \(f_i,f_j\) multiply as expected within the ring, \(\mathbb{Z}[q,q^{-1}]\). Moreover, the braids \(B_i\) and \(Y_j\) concatenate as shown above. 

\begin{defn}

There exists a function know as the Jones-Ocneanu trace, \(Tr: \mathbb{H}_n \to \mathbb{Q}(q,a)\), that is evaluated by:

1. \[ \sum_{i=0}^{m} f_i \vcenteredinclude{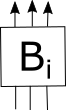} \to \sum_{i=0}^{m} f_i \vcenteredinclude{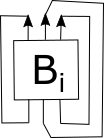}, f_i \in \mathbb{Z}[q,q^{-1}], B_i \in B_n \]

2. \[\vcenteredincludeb{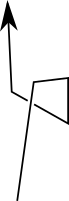}=-q^{-1}\vcenteredincludeb{loud4}, \vcenteredincludeb{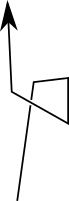}=-q^{-1}a^2\vcenteredincludeb{loud4}, \vcenteredincludeb{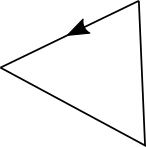}= \frac{1-a^2}{1-q^2}\]

\end{defn}

Definition 2.16 describes the Jones-Ocneanu trace function which maps an element in the Hecke algebra to an element in the field of rational functions in q and a. This is first done by mapping a  Hecke algebra element to a \(\mathbb{Z}[q,q^{-1}]\)-module through the closure of the braids in the Hecke algebra element seen in step 1. Next, step 2 allows one to reduce this new element to an element in \(\mathbb{Q}(q,a)\). Moreover, its important to notice that the Hecke algebra relations seen in Definition 2.15 will also be important in the Jones-Ocneanu trace evaluation.

\begin{ex}

Calculating the trace of \(\vcenteredincludeb{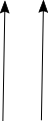} \in \mathbb{H}_2\):
\\
1. \(\vcenteredincludeb{ex1} \to \vcenteredincludeb{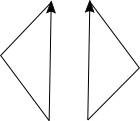}\)
\\
2. \(\vcenteredincludeb{ex2} = (\frac{1}{q-q^{-1}})(\vcenteredincludeb{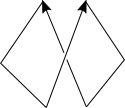} - \vcenteredincludeb{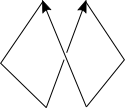})\)
\\
3. \(\vcenteredincludeb{ex2} = (\frac{1}{q-q^{-1}})((-q^{-1})\vcenteredincludeb{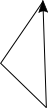} - (-q^{-1}a^2)\vcenteredincludeb{ex5})\)
\\
4. \(\vcenteredincludeb{ex2} = (\frac{1}{q-q^{-1}})((-q^{-1})(\frac{1-a^2}{1-q^2}) - (-q^{-1}a^2)(\frac{1-a^2}{1-q^2}))\)

\end{ex}

Example 2.17 walks through the trace calculation of a simple braid element in \(\mathbb{H}_2\). Step 1 shows the closure of the element. Step 2 utilizes the Hecke algebra relations to reduce the no crossing to a positive crossing and a negative crossing. Lastly, steps 3 and 4 utilize the trace evaluation relations in order to reduce the element down further. 

\subsubsection{HOMFLY Polynomial}

\begin{defn}

Similar to the Jones polynomial, the HOMFLY polynomial is another link invariant that can be calculated utilizing the Jones-Ocneanu trace in the following way:

1. Write \(L\) as the closure of a braid \(B\) (Alexander's Theorem).

2. Calculate the \(Tr(B)\) for \(B \in \mathbb{H}_n\)

3. \(HOMFLY(L) = (qa^{-1})^{n(B)}(-a)^{-w(B)}Tr(B)\)

\(n(B)\) = number of strands in \(B\)

\(w(B) = \) number of positive crossings in \(B\) minus the number of negative crossings in \(B\)

\end{defn}

In the next section, I will utilize the Jones-Ocneanu trace in order to calculate a general formula for the HOMFLY polynomial of the closure of a specific family of braids.

\section{Independent Research}

In each braid group there is one coxeter braid which permutes the left most node to the right most node through positive crossings. Each other node is subsequently permuted to the left by one node. The coxeter braid of rank 4 can be seen in Figure 15, below.

\begin{figure}[h]
    \centering
    \includegraphics[width=0.4\textwidth]{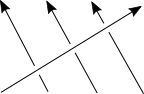}
    \caption{Coxeter Braid in \(B_4\)}
\end{figure}

\begin{figure}[h]
    \centering
    \includegraphics[width=0.4\textwidth]{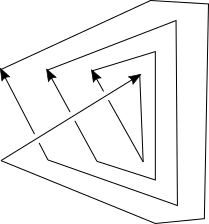}
    \caption{Closure of Coxeter Braid in \(B_4\)}
\end{figure}

The closure of coxeter braids, seen in Figure 16, is nothing more than the unknot. This can be shown by using multiple type 1 Reidemeister moves in order to undo the positive crossings in the braid closure. A braid that is a bit more interesting than the coxeter braid, which will be referred to as the looped coxeter braid, can be seen in Figure 17.

\begin{figure}[h]
    \centering
    \includegraphics[width=0.4\textwidth]{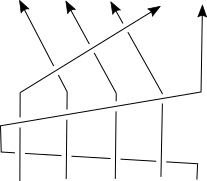}
    \caption{Looped Coxeter Braid in \(B_5\)}
\end{figure}

Each coxeter braid of rank n corresponds to a single looped coxeter braid of rank n+1. Therefore, there is no rank 1 looped coxeter braid. An example of the closure of these looped coxeter braids can be seen in Figure 18.

\begin{figure}[h]
    \centering
    \includegraphics[width=0.25\textwidth]{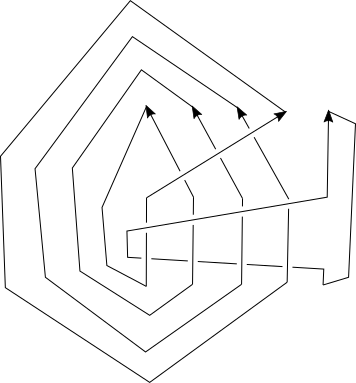}
    \caption{Looped Coxeter Braid in \(B_5\) Closure}
\end{figure}

As shown, the closure produces a link with two components. Those components being an unknot that raps around a coxeter braid closure. In order to develop a general equation for the closure of looped coxeter braids, one needs to first inspect the Jones-Ocneanu trace evaluations for these braids. This is shown in the next section.

\subsection{Looped Coxeter Braid Trace Evaluation}

First, it is important to manipulate a general looped coxeter braid within its corresponding Hecke algebra. The figures used will visually show a looped coxeter braid of rank 5, but this is just used as a place holder for a looped coxeter braid of rank n, where \(n > 3\). The Hecke algebra looped coxeter element is manipulated as follows:

1. \[\vcenteredincludeb{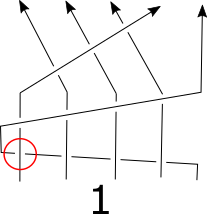} = (q-q^{-1})\vcenteredincludeb{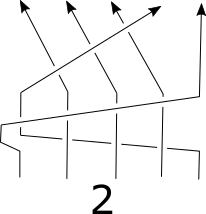} + \vcenteredincludeb{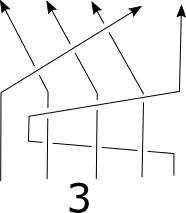}\]

2. \[(q-q^{-1})\vcenteredincludeb{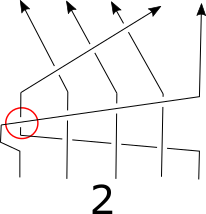} = (q-q^{-1})^2\vcenteredincludeb{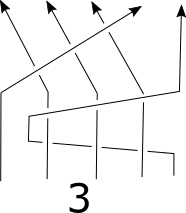} + (q-q^{-1})\vcenteredincludeb{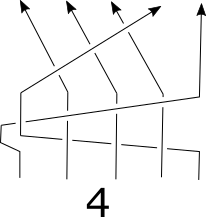}\]

3. \[(q-q^{-1})\vcenteredincludeb{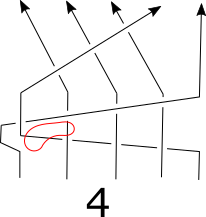} = \vcenteredincludeb{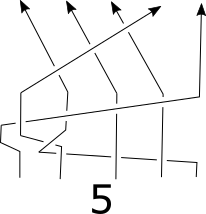} - \vcenteredincludeb{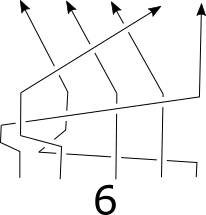}\]

4. \[\vcenteredincludeb{coxeter3} = ((q-q^{-1})^2+1)(Braid 3)+(Braid 5)-(Braid 6)\]

In the first line, we use the first Hecke algebra relation in order to equate our looped coxeter braid of rank n to braids 2 and 3. Through this relation, the positive crossing outlined in red is changed to a negative crossing in braid 3 and a no crossing in braid 2. In the second line, we use this same relation to change the positive crossing outlined in red in braid 2 to a no crossing, which gives us braid 3 again, and a negative crossing seen in braid 4. In the third line, we equate braid 4 to braids 5 and 6 by changing the no crossing outlined in red to a positive crossing and a negative crossing. Lastly, the fourth line displays how we can write our general looped coxeter braid as a linear combination of braids 3, 5, and 6. It will be useful to use this equal element when computing the trace evaluation of the general looped coxeter braid. We can do this by first systematically evaluating the trace evaluation of the closure of braids 3,5, and 6. Beginning with braid 3, we get that

\[\vcenteredincluded{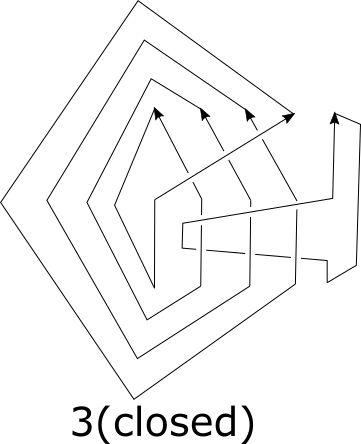} = (-q^{-1})\vcenteredincluded{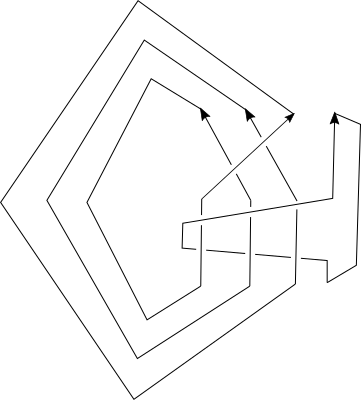}.\]

We see that the link on the right is actually the closure of the looped coxeter braid of rank n-1. Therefore, we get that the

\[Tr(Braid 3) = (-q^{-1})Tr(LCB_{n-1}).\]

Here, 'LCB' stands for looped coxeter braid and the subscript denotes the rank. Moving on to braid 5, it can be seen that

\[\vcenteredincluded{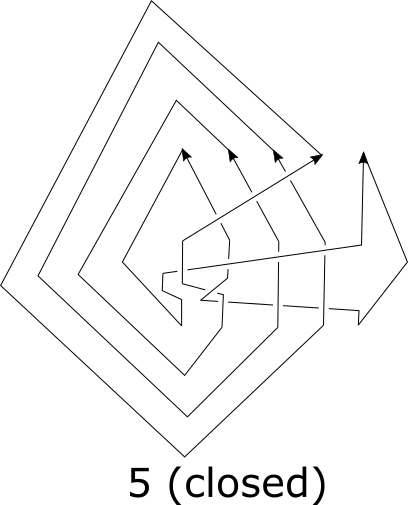} = (-q^{-1})\vcenteredincluded{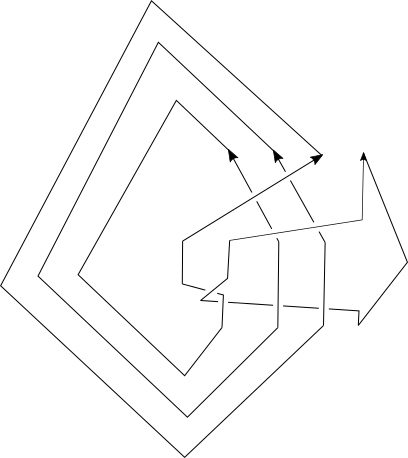}.\]

Similarly to braid 3, we get that the

\[Tr(Braid 5) = (-q^{-1})Tr(LCB_{n-1}).\]

Lastly, in regards to braid 6, we have that

\[\vcenteredincluded{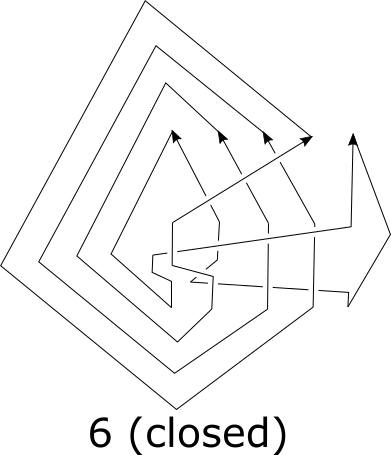} = \vcenteredincluded{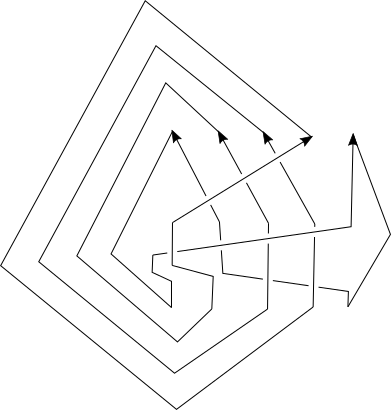}.\]

Simplifying further by using our type 1 Reidemeister move relation, we get that

\[\vcenteredincluded{coxeter21} = (-q^{-1})\vcenteredincluded{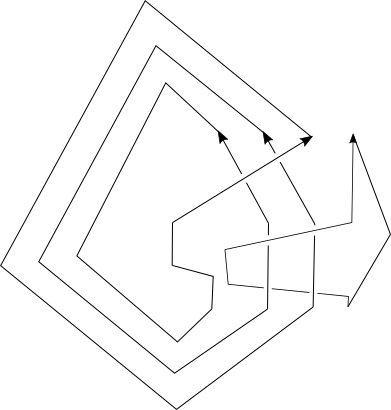}\]

and

\[(-q^{-1})\vcenteredincluded{coxeter22} = (q^{-2})\vcenteredincluded{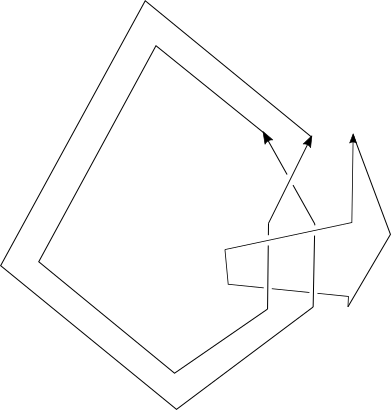}.\]

Notice that the final link on the bottom right is actually the braid closure of the looped coxeter braid of rank n-2. Therefore, we can write that the

\[Tr(Braid 6) = (q^{-2})Tr(LCB_{n-2}).\]

Putting all of the trace evaluations together, we get that the

\[Tr(LCB_{n}) = (-q-q^{-3})Tr(LCB_{n-1})-(q^{-2})Tr(LCB_{n-2}).\]

As seen, we end up with a recursive general formula for the trace evaluation of the looped coxeter braid. This can be combined with the number of strands in the looped coxeter braid and the number of positive crossings minus the number of negative crossings to form a general formula for the HOMFLY polynomial of the closed looped coxeter braid:

\[HOMFLY(CLCB_n) = (qa^{-1})^n(-a)^{-3n+4}Tr(LCB_n)\]

\section{Plan for Future Work}

The Jones-Ocneanu trace is not the only trace function that exists. In fact, a function \(g: \mathbb{H}_n \to V\), where \(V\) is some  \(\mathbb{Z}[q,q^{-1}]-module\), is considered a trace function if \(g(\alpha\beta) = g(\beta\alpha)\) for the braids \(\beta,\alpha \in \mathbb{H}_n\) and if the function is \(\mathbb{Z}[q,q^{-1}]-linear\). This means that

\[g(\sum_{i=0}^{m} f_i \vcenteredinclude{braidb}) = \sum_{i=0}^{m} f_i g(\vcenteredinclude{braidb}).\]

A special trace function, referred to as the universal trace, is evaluated by mapping the braid elements in \(\mathbb{H}_n\) to the closure of these elements in the punctured plane. This can be pictured as follows:

\[\sum_{i=0}^{m} f_i \vcenteredinclude{braidb} \to \sum_{i=0}^{m} f_i \vcenteredinclude{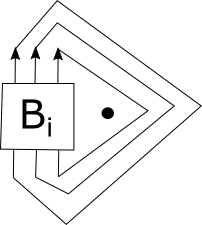}\]

It is important that each corresponding strand is closed around the punctured point. Moreovr, the \(\mathbb{Z}[q,q^{-1}]-module\) that the universal trace maps to is referred to as \(\Lambda_n\). The reason why the universal trace is referred to as 'universal' is because all other trace functions can be written as first mapping to \(\Lambda_n\) through the universal trace, and then mapping from \(\Lambda_n\) to the desired \(\mathbb{Z}[q,q^{-1}]-module\). Due to this feature, having a deep understanding of the universal trace aids in understanding other trace functions like the Jones-Ocneanu trace. My research also becomes relevant due to the fact that all braid closures in \(\Lambda_n\) can be written as the closure of a "linear combination" of coxeter braids (a rank 6 braid that is a linear combination of a coxeter braid of rank 3, 2, and 1 is shown in Figure 17).

\begin{figure}[h]
    \centering
    \includegraphics[width=0.4\textwidth]{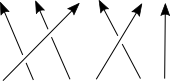}
    \caption{Closure of Coxeter Braid in \(B_4\)}
\end{figure}

Moreover, there are many important operators defined on \(\Lambda_n\) that involve adding loops similar to the looped coxeter braid that I looked at in my research. In conclusion, further pursuit of my research would involve studying the universal trace and the effect that these certain operators have on elements in \(\Lambda_n\).
\\
\\
Note: All unreferenced images in Literature Review, Independent Research, and Plan for Future Work was created by me using Inkspace.

\section{Acknowledgement}

I would like to give thanks to Professor  Anthony Iarrobino for aiding me in the mathematics, writings, and presentations required for this project. I would also like to thank my consultant, Professor Matt Hogancamp, for taking the time to help me navigate some of the mathematics involved in this project.




\end{document}